\documentclass{amsart}

\usepackage{fancyhdr}
\usepackage{lastpage}
\usepackage{stmaryrd,yhmath}
\usepackage{tikz}

\newcommand{\bdis}{\begin{displaymath}}
\newcommand{\edis}{\end{displaymath}}
\newcommand{\be}{\begin{equation}}
\newcommand{\ee}{\end{equation}}
\newcommand{\mbb}{\mathbb}
\newcommand{\mcal}{\mathcal}

\theoremstyle{definition}

\theoremstyle{remark}
\newtheorem{remark}[]{Remark}

\newtheorem*{mydef11}{{\bf Theorem 1}}

\newtheorem*{mydef12}{{\bf Theorem 2}}

\newtheorem*{mydef2}{{\bf Definition}}

\newtheorem*{mydef4}{{\bf Corollary}}

\newtheorem*{mydef51}{{\bf Lemma 1}}

\newtheorem*{mydef52}{{\bf Lemma 2}}

\numberwithin{equation}{section}

\begin{document}

\title{Jacob's ladders and some unbounded decomposition of the $\zeta$-oscillating systems on products of other $\zeta$-oscillating systems as
an $\zeta$-analogue of the prime decomposition}

\author{Jan Moser}

\address{Department of Mathematical Analysis and Numerical Mathematics, Comenius University, Mlynska Dolina M105, 842 48 Bratislava, SLOVAKIA}

\email{jan.mozer@fmph.uniba.sk}


\begin{abstract}
In this paper we introduce new class of multiplicative interactions of the $\zeta$-oscillating systems generated by a subset of power functions. The main result
obtained expresses an analogue of prime decomposition (without the property of uniqueness).  \\

Dedicated to recalling of Nicola Tesla's oscillators
\end{abstract}
\maketitle

\section{Introduction}

\subsection{}

Let us remind that in our papers \cite{1} -- \cite{8} we have introduced within the theory of the Riemann zeta-function the following notions: Jacob ladders (JL),
$\zeta$-oscillating systems (OS), factorization formula (FF), metamorphosis of the oscillating systems (M), $\mcal{Z}_{\zeta,Q^2}$-transformations (ZT), and
interactions between oscillating systems (IOS).

The main result obtained in this direction (IOS) in the paper \cite{8} is the following set of $\zeta$-analogues of the elementary trigonometric identity $\cos^2 t+\sin^2 t=1$:
\be \label{1.1}
\begin{split}
& \cos^2(\alpha_0^{2,2})\prod_{r=1}^{k_2}\left|\frac{\zeta\left(\frac 12+i\alpha_r^{2,2}\right)}{\zeta\left(\frac 12+i\beta_r^2\right)}\right|^2+\\
& +
\sin^2(\alpha_0^{1,1})\prod_{r=1}^{k_1}\left|\frac{\zeta\left(\frac 12+i\alpha_r^{1,1}\right)}{\zeta\left(\frac 12+i\beta_r^1\right)}\right|^2\sim 1,\ L\to\infty,
\end{split}
\ee
(for the notations used above see \cite{8}).

\begin{remark}
We may call the kind of interactions of the $\zeta$-oscillating systems in (\ref{1.1}) as the linear interaction between the systems (LIOS), comp. section 2.4 in \cite{8}.
\end{remark}

Of course, in our paper \cite{8}, we have obtained also other formulae for interactions of $\zeta$-oscillating systems that are also good characterized by the words
\emph{linear interactions}, as in Remark 1.

\subsection{}

In this paper we obtain different kind of the formula w.r.t. (\ref{1.1}). Namely, we obtain new formula that is essentially non-linear multiplicative in corresponding $\zeta$-oscillating systems:
\be \label{1.2}
\begin{split}
& \prod_{r=1}^k\left|\frac{\zeta\left(\frac 12+i\alpha_r\right)}{\zeta\left(\frac 12+i\beta_r\right)}\right|^2\sim \\
& \sim \frac{1}{\Delta+1}\prod_{l=1}^n (\Delta_l+1)\frac{1}{(\alpha_0-L)^\Delta}\prod_{l=1}^n(\alpha_0^{\Delta_l,k_l}-L)^{\Delta_l}\times \\
& \times \prod_{l=1}^n\prod_{r=1}^{k_l}\left|\frac{\zeta\left(\frac 12+i\alpha_r^{\Delta_l,k_l}\right)}{\zeta\left(\frac 12+i\beta_r^{k_l}\right)}\right|^2,\
L\in\mbb{N},\ L\to\infty,
\end{split}
\ee
where
\bdis
1\leq k_l\leq k_0,\ \Delta=\sum_{l=1}^n \Delta_l,\ \Delta>\Delta_1\geq \Delta_2\geq \dots \geq \Delta_n,
\edis
\bdis
\Delta,\Delta_l\in\mbb{R}^+,\ k_0,n\in\mbb{N},
\edis
for every fixed $\Delta,k_0,n$.

\begin{remark}
A new property of interactions between the elements of set of $\zeta$-oscillating systems generated over subset of class of power functions is expressed by our formula (\ref{1.2}). Namely,
given the main (say) $\zeta$-oscillating system
\be \label{1.3}
\prod_{r=1}^k\left|\frac{\zeta\left(\frac 12+i\alpha_r\right)}{\zeta\left(\frac 12+i\beta_r\right)}\right|^2\longleftarrow \Delta
\ee
is now factorized itself (comp. the factorization formula (3.8) in \cite{8}) by means of the choice of basic system of the $\zeta$-oscillating systems
\be \label{1.4}
\prod_{r=1}^{k_l}\left|\frac{\zeta\left(\frac 12+i\alpha_r^{\Delta_l,k_l}\right)}{\zeta\left(\frac 12+i\beta_r^{k_l}\right)}\right|^2\longleftarrow \Delta_l,\ l=1,\dots,n.
\ee
\end{remark}

\begin{remark}
Let us remind the prime decomposition
\be \label{1.5}
n=p_1\cdots p_l,\ n\in\mbb{N}.
\ee
Our formula (\ref{1.2}) represents, in this direction, decomposition of the main $\zeta$-oscillating system (\ref{1.3}) into weighted product of basic $\zeta$-oscillating systems
(\ref{1.4}). That is the formula (\ref{1.2}) represents a $\zeta$-analogue of the prime decomposition (\ref{1.5}), except its uniqueness.
\end{remark}

Let us remind the diagram (see \cite{8}) of corresponding notions we have introduced within the theory of the Riemann zeta-function. The last lot in this diagram is IOS = interaction between
oscillating systems. Now, the mentioned lot itself may be split into two others as follows
\bdis
\mbox{IOS}\rightarrow \left\{\begin{array}{l} \mbox{LIOS} \\ \mbox{NIOS},   \end{array} \right.
\edis
about LIOS see Remark 1 and NIOS means non-linear IOS (for example, (\ref{1.2})).

\section{Lemma}

\subsection{}

This paper is devoted mainly to study of some type of multiplicative interactions between the $\zeta$-oscillating systems generated by the following subset of the class of power functions
\be \label{2.1}
\begin{split}
& f(t)=f(t;L,\Delta)=(t-L)^\Delta, \\
& t\in [L,L+U],\ l\in\mbb{N},\ \Delta>0, \\
& U\in (0,a],\ a\in (0,1).
\end{split}
\ee

\begin{remark}
The assumption
\bdis
U\in (0,a], a<1
\edis
in (\ref{2.1}) was chosen because of the possibility of interpretation of our results in terms of deterministic signals (pulses) in the communication theory. In this connection, see also
the notion of the $\mcal{Z}_{\zeta,Q^2}$-transformation (device) we have introduced in our paper \cite{7}.
\end{remark}

Since
\bdis
f(L;L,\Delta)=0
\edis
then we use, instead of Definition 2 from \cite{8}, the following

\begin{mydef2}
The symbol
\bdis
f(t)\in \tilde{C}_0[T,T+U]
\edis
stands for the following
\be \label{2.2}
\begin{split}
& f(t)\in C[T,T+U] \ \wedge \ f(t)\geq 0 \ \wedge \\
& \wedge \ \exists\- t_0\in [T,T+U]:\ f(t_0)>0, \\
& T>T_0, \ U\in (0,U_0],\ U_0=o\left( \frac{T}{\ln T}\right),\ T\to\infty.
\end{split}
\ee
\end{mydef2}

\begin{remark}
Of course,
\bdis
\tilde{C}\subset \tilde{C}_0,
\edis
and
\bdis
f(t;L,\Delta)\in\tilde{C}_0[L,L+U],\ L>T_0,\ \Delta>0.
\edis
\end{remark}

\begin{remark}
The case
\bdis
\Delta\in (-1,0)
\edis
is excluded though
\bdis
\int_L^{L+U} (t-L)^\Delta{\rm d}t<+\infty,
\edis
because
\bdis
(t-L)^\Delta\rightarrow+\infty \ \mbox{as} \ t\to L^+ ,
\edis
(the case $\Delta=0$ is the trivial one).
\end{remark}

\begin{remark}
We see that (\ref{2.2}) implies the following
\bdis
\int_T^{T+U} f(t){\rm d}t>0,
\edis
and this is sufficient for applicability of our algorithm to generate of factorization formula (comp. \cite{8}, (3.1) -- (3.11)).
\end{remark}

\subsection{}

Since (see (\ref{2.1}))
\be \label{2.3}
\frac 1U\int_L^{L+U} (t-L)^\Delta{\rm d}t=\frac{1}{\Delta+1}U^\Delta,\ \Delta>0 ,
\ee
then we have by just mentioned algorithm (see Remark 7) the following

\begin{mydef51}
For the function (\ref{2.1}) there are vector-valued functions
\bdis
(\alpha_0,\alpha_1,\dots,\alpha_k,\beta_1,\dots,\beta_k),\ k=1,\dots,k_0,\ k_0\in\mbb{N},
\edis
($k_0$ being arbitrary and fixed) such that the following factorization formula holds true
\be \label{2.4}
\begin{split}
& \prod_{r=1}^k\left|\frac{\zeta\left(\frac 12+i\alpha_r\right)}{\zeta\left(\frac 12+i\beta_r\right)}\right|^2\sim
\frac{1}{\Delta+1}\left(\frac{U}{\alpha_0-L}\right)^\Delta, \\
& L\in\mbb{N}, \ L\to\infty,\ U\in (0,a], \ a<1,\ \Delta>0,
\end{split}
\ee
where
\be \label{2.5}
\begin{split}
& \alpha_r=\alpha_r(U,L,\Delta,k),\ r=0,1,\dots,k, \\
& \beta_r=\beta_r(U,L,k),\ r=1,\dots,k, \\
& L<\alpha_0<L+U \ \Rightarrow \ 0<\alpha_0-L<U.
\end{split}
\ee
\end{mydef51}

\begin{remark}
Of course, in the asymptotic formula (\ref{2.4}) the symbol $\sim$ stands for (comp. \cite{8}, (3.8))
\bdis
=\left\{ 1+\mcal{O}\left(\frac{\ln\ln L}{\ln L}\right)\right\}(\dots).
\edis
\end{remark}

\section{On unbounded decomposition of $\zeta$-oscillating systems}

\subsection{}

Now, we have, by our Lemma 1, that
\be \label{3.1}
\begin{split}
& (\Delta_l+1)(\alpha_0^{\Delta_l,k_l}-L)\prod_{r=1}^{k_l}
\left|\frac{\zeta\left(\frac 12+i\alpha_r^{\Delta_l,k_l}\right)}{\zeta\left(\frac 12+i\beta_r^{k_l}\right)}\right|^2\sim U^{\Delta_l}, \\
& L\to\infty, \ \Delta_l>0,\ 1\leq k_l\leq k_0,
\end{split}
\ee
where (see (\ref{2.5}))
\bdis
\Delta\rightarrow\Delta_l,\ l\rightarrow k_l,\ \Rightarrow \ \alpha_r\rightarrow \alpha_r^{\Delta_l,k_l}, \dots
\edis
Consequently, we obtain by Lemma 1 and (\ref{3.1}) the following

\begin{mydef11}
In the case
\be \label{3.2}
\begin{split}
& \Delta=\sum_{l=1}^n \Delta_l,\ \Delta>\Delta_1\geq \dots \geq \Delta_n, \\
& \Delta,\Delta_l>0,\ n\in\mbb{N} ,
\end{split}
\ee
there are vector-valued functions
\bdis
\begin{split}
& (\alpha_0,\alpha_1,\dots,\alpha_k,\beta_1,\dots,\beta_k), \\
& (\alpha_0^{\Delta_l,k_l},\alpha_1^{\Delta_l,k_l},\dots,\alpha_k^{\Delta_l,k_l},\beta_1^{k_l},\dots,\beta_k^{k_l}), \\
& l=1,\dots,n
\end{split}
\edis
where
\bdis
\begin{split}
& \alpha_r=\alpha_r(U,L,\Delta,k),\ r=0,1,\dots,k, \\
& \beta_r=\beta_r(U,L,k),\ r=1,\dots,k, \\
& \alpha_r^{\Delta_l,k_l}=\alpha_r^{\Delta_l,k_l}(U,L,\Delta_l,k_l),\ r=0,1,\dots,k_l, \\
& \beta_r^{k_l}=\beta_r^{k_l}(U,L,k_l),\ r=1,\dots,k_l,
\end{split}
\edis
such that
\be \label{3.3}
\begin{split}
& \prod_{r=1}^k\left|\frac{\zeta\left(\frac 12+i\alpha_r\right)}{\zeta\left(\frac 12+i\beta_r\right)}\right|^2\sim \\
& \sim \frac{1}{\Delta+1}\prod_{l=1}^n(\Delta_l+1)\frac{1}{(\alpha_0-L)^\Delta}\prod_{l=1}^n(\alpha_0^{\Delta_l,k_l}-L)^{\Delta_l}\times \\
& \times \prod_{l=1}^n\prod_{r=1}^{k_l}
\left|\frac{\zeta\left(\frac 12+i\alpha_r^{\Delta_l,k_l}\right)}{\zeta\left(\frac 12+i\beta_r^{k_l}\right)}\right|^2, \\
& L\in\mbb{N},\ L\to\infty,
\end{split}
\ee
where
\bdis
\begin{split}
& L<\alpha_0<L+U \ \Rightarrow \ 0<\alpha_0-L<U, \\
& L<\alpha_0^{\Delta_l,k_l}<L+U \ \Rightarrow \ 0<\alpha_0^{\Delta_l,k_l}-L<U,
\end{split}
\edis
for arbitrary but fixed $k_0,n$.
\end{mydef11}

\begin{remark}
The assumption of arbitrary but fixed $n$ has its reason in the following
\bdis
\prod_{l=1}^n \{ 1+o(1)\}=1+o(1),\ L\to\infty
\edis
that must be obeyed.
\end{remark}

\subsection{}

First, we give the following

\begin{remark}
In our paper \cite{8} we have introduced some classes of additive interactions (linear ones in the corresponding non-linear $\zeta$-oscillating
systems). In contrast with these, our formula (\ref{3.3}) expresses some class of multiplicative interactions of the corresponding $\zeta$-oscillating
systems.
\end{remark}

Secondly, the formula (\ref{3.3}) contains following products:

\begin{itemize}
\item[(a)] the main (say) $\zeta$-oscillating system
\be \label{3.4}
\prod_{r=1}^k\left|\frac{\zeta\left(\frac 12+i\alpha_r\right)}{\zeta\left(\frac 12+i\beta_r\right)}\right|^2 \leftarrow \Delta\in\mbb{R}^+,
\ee
\item[(b)] the basic set (say) of the $\zeta$-oscillating systems
\be \label{3.5}
\left|\frac{\zeta\left(\frac 12+i\alpha_r^{\Delta_l,k_l}\right)}{\zeta\left(\frac 12+i\beta_r^{k_l}\right)}\right|^2 \leftarrow \Delta_l\in\mbb{R}^+,\ l=1,\dots , n ,
\ee
\item[(c)] the product
\be \label{3.6}
\frac{1}{\Delta+1}\prod_{l=1}^n (\Delta_l+1)
\ee
that contains the generating parameters (see (\ref{3.2})),
\item[(d)] the product
\be \label{3.7}
\frac{1}{(\alpha_0-L)^\Delta}\prod_{l=1}^n(\alpha_0^{\Delta_l,k_l}-L)^{\Delta_l}
\ee
that contains the control functions
\be \label{3.8}
\alpha_0,\ \alpha_0^{\Delta_l,k_l}.
\ee
\end{itemize}

\begin{remark}
The notion of \emph{control functions} (\ref{3.8}) is based on the fact that the quotient of the $\zeta$-oscillating system (\ref{3.4}) with the product
\bdis
\prod_{l=1}^n\prod_{r=1}^{k_l}
\left|\frac{\zeta\left(\frac 12+i\alpha_r^{\Delta_l,k_l}\right)}{\zeta\left(\frac 12+i\beta_r^{k_l}\right)}\right|^2
\edis
is proportional to the factor (\ref{3.7}); the product (\ref{3.6}) is constant for fixed representation (\ref{3.2}) that is for corresponding interaction.
\end{remark}

\begin{remark}
Consequently, the formula (\ref{3.3}) expresses a kind of superfactorization of the main $\zeta$-oscillating system (\ref{3.4}). Namely, this main system is factorized
itself by means of the basic of others $\zeta$=oscillating systems (\ref{3.5}). Of course, it is true that all the oscillating systems playing role are generated by the class
of used power functions with $\Delta,\Delta_l\in\mbb{R}^+$.
\end{remark}

\subsection{}

Now, we give, for example, the following extremal case of the superfactorization (\ref{3.3}).

\begin{mydef4}
If
\bdis
k=k_1=k_2=\dots=k_n=1,
\edis
then
\be \label{3.9}
\begin{split}
& \left|\frac{\zeta\left(\frac 12+i\bar{\alpha}_1\right)}{\zeta\left(\frac 12+i\beta_1\right)}\right|^2\sim  \\
& \sim \frac{1}{\Delta+1}\prod_{l=1}^n(\Delta_l+1)\frac{1}{(\bar{\alpha}_0-L)^\Delta}\prod_{l=1}^n (\alpha_0^{\Delta_l,1}-L)^{\Delta_l}\times \\
& \times \prod_{l=1}^n \left|\frac{\zeta\left(\frac 12+i\alpha_1^{\Delta_l,1}\right)}{\zeta\left(\frac 12+i\beta_1^{1}\right)}\right|^2,\ L\to\infty,
\end{split}
\ee
where
\bdis
\begin{split}
& \bar{\alpha}_r=\alpha_r(U,L,\Delta,1),\ r=0,1, \\
& \beta_1=\beta_1^1=\beta_1(U,L,1).
\end{split}
\edis
\end{mydef4}

\begin{remark}
Other extremal formulae correspond to the cases:
\bdis
\begin{split}
& k=k_0,\ k_1=k_2=\dots=k_n=1, \\
& k=1,\ k_1=k_2=\dots=k_n=k_0, \\
& k=k_1=\dots=k_n=k_0.
\end{split}
\edis
\end{remark}

\section{Some additive interactions}

\subsection{}

If
\be \label{4.1}
\begin{split}
& f(t)=f(t;L,\Delta_1,\dots,\Delta_n)=\sum_{l=1}^n (t-L)^{\Delta_l}, \\
& t\in [L,L+U],\ \Delta_l>0,
\end{split}
\ee
of course,
\bdis
f(t)\in \tilde{C}_0[L,L+U]
\edis
(comp. Definition), then
\bdis
\frac 1U\int_L^{L+U}\sum_{l=1}^n (t-L)^{\Delta_l} {\rm d}t=\sum_{l=1}^n\frac{1}{\Delta_l+1}U^{\Delta_l},
\edis
and, consequently, we obtain by our algorithm the following

\begin{mydef52}
For the function (\ref{4.1}) there are vector-valued functions
\bdis
(\tilde{\alpha}_0,\tilde{\alpha}_1,\dots,\tilde{\alpha}_k,\tilde{\beta}_1,\dots,\tilde{\beta}_k),\ k=1,\dots,k_0
\edis
such that the following factorization formula
\be \label{4.2}
\prod_{r=1}^k\left|\frac{\zeta\left(\frac 12+i\tilde{\alpha}_r\right)}{\zeta\left(\frac 12+i\tilde{\beta}_r\right)}\right|^2\sim
\frac{\sum_{l=1}^n\frac{1}{\Delta_l+1}U^{\Delta_l}}{\sum_{l=1}^n (\tilde{\alpha}_0-L)^{\Delta_l}},\ \Delta_l>0,\ L\to\infty
\ee
holds true, where
\bdis
\begin{split}
& \tilde{\alpha}_r=\alpha_r(U,L,\Delta_1,\dots,\Delta_n,k),\ r=0,1,\dots,k, \\
& \tilde{\beta}_r=\beta_r(U,L,1),\ r=1,\dots,k..
\end{split}
\edis
\end{mydef52}

Next, we have (see (\ref{3.1})) the following formula
\be \label{4.3}
\begin{split}
& (\alpha_0^{\Delta_l,k_l}-L)^{\Delta_l}\prod_{r=1}^{k_l}
\left|\frac{\zeta\left(\frac 12+i\alpha_r^{\Delta_l,k_l}\right)}{\zeta\left(\frac 12+i\beta_r^{k_l}\right)}\right|^2\sim \frac{1}{\Delta_l+1}U^{\Delta_l}, \\
& \Delta_l>0,\ 1\leq k_l\leq k_0.
\end{split}
\ee
Consequently, we have (see (\ref{4.2}), (\ref{4.3})) the following

\begin{mydef12}
\be \label{4.4}
\begin{split}
& \prod_{r=1}^k\left|\frac{\zeta\left(\frac 12+i\tilde{\alpha}_r\right)}{\zeta\left(\frac 12+i\tilde{\beta}_r\right)}\right|^2\sim \\
& \sim \frac{\sum_{l=1}^n (\alpha_0^{\Delta_l,k_l}-L)^{\Delta_l}}{\sum_{l=1}^n (\tilde{\alpha}_0-L)^{\Delta_l}}
\prod_{r=1}^{k_l}
\left|\frac{\zeta\left(\frac 12+i\alpha_r^{\Delta_l,k_l}\right)}{\zeta\left(\frac 12+i\beta_r^{k_l}\right)}\right|^2,\ L\to\infty, \\
& \Delta_l>0,\ l=1,\dots,n.
\end{split}
\ee
\end{mydef12}

\begin{remark}
The additive interaction (comp. \cite{8}) between $\zeta$-oscillating systems
\bdis
\prod_{r=1}^k\left|\frac{\zeta\left(\frac 12+i\tilde{\alpha}_r\right)}{\zeta\left(\frac 12+i\tilde{\beta}_r\right)}\right|^2 \leftarrow \{ \Delta_1,\dots,\Delta_n\},
\edis
is expressed by the formula (\ref{4.4}), and the set of $n$ others oscillating systems
\bdis
\prod_{r=1}^{k_l}
\left|\frac{\zeta\left(\frac 12+i\alpha_r^{\Delta_l,k_l}\right)}{\zeta\left(\frac 12+i\beta_r^{k_l}\right)}\right|^2 \leftarrow \Delta_l,\ l=1,\dots,n.
\edis
Of course, we also may write the set of $n$ similar formulae to (\ref{4.4}).
\end{remark}

\begin{remark}
Let us notice explicitly that the assumption (see (\ref{4.4}))
\bdis
\Delta_l>0,\ l=1,\dots,n
\edis
on the set
\bdis
\{ \Delta_1,\dots,\Delta_n\}
\edis
is the only one (in comparison with the Theorem 1).
\end{remark}


\end{document}